\documentclass{amsart}
\usepackage{latexsym,a4,mathrsfs,amsthm,amsmath,amssymb}
\usepackage[latin1]{inputenc}
\usepackage{graphics}
\usepackage{epsfig}
\usepackage[english]{babel}
\usepackage{hyperref}
\usepackage[all]{xy}
\usepackage{verbatim}
\usepackage[dvipdfm]{geometry}
\usepackage{fancyhdr}
\usepackage{lastpage}
\newtheorem{theorem}{Theorem}[section]
\newtheorem{lemma}[theorem]{Lemma}

\newtheorem{definition}[theorem]{Definition}

\newtheorem{corollary}{Corollary}[theorem]
\theoremstyle{definition}
\newtheorem{example}[theorem]{Example}
\numberwithin{equation}{section}

\newcommand{\rw}{\rightarrow}

\newcommand{\Rr}{\mathbb{R}}
\newcommand{\C}{\mathbb{C}}
\newcommand{\Q}{\mathbb{Q}}
\newcommand{\Zr}{\mathbb{Z}}
\newcommand{\Nr}{\mathbb{N}}
\newcommand{\eps}{\epsilon}
\newcommand{\la}{\lambda}
\newcommand{\Log}{\mathrm{Log}}

\newcommand{\bx}{\mathbf{x}}
\newcommand{\bn}{\mathbf{n}}
\newcommand{\bs}{\mathbf{s}}

\newcommand{\bk}{\mathbf{k}}
\newcommand{\norm}[1]{\| #1 \|}
\renewcommand{\Re}{\mathrm{Re}}
\renewcommand{\Im}{\mathrm{Im}}
\newcommand{\nth}[1]{\left[ #1 \right]_n}
\newcommand{\ud}[1]{\underline{#1}}

\newcommand{\set}[1]{\left\{#1\right\}}
\newcommand{\E}[1]{\mathbb{E}\left[#1\right]}
\newcommand{\Pb}[1]{\mathbb{P}\left[#1\right]}
\begin{document}


\title[Permutation matrices and their characteristic polynomials]{Permutation matrices and the moments of their characteristics polynomials}




\address{%
Department of Mathematics, University Z\"urich, Winterthurerstrasse 190, 8057 Z\"urich, CH
}
\email{dirkz@math.uzh.ch}

\begin{abstract}
In this paper, we are interested in the moments of the characteristic polynomial $Z_n(x)$ of the $n\times n$ permutation matrices with respect to the uniform measure.
We use a combinatorial argument to write down the generating function of $\E{\prod_{k=1}^p Z_n^{s_k}(x_k)}$ for $s_k\in\Nr$. We show with this generating function that $\lim_{n\rw\infty} \E{\prod_{k=1}^p Z_n^{s_k}(x_k)}$ exists for $\max_k|x_k|<1$ and calculate the growth rate for $p=2, |x_1|=|x_2|=1$, $x_1=\overline{x_2}$ and $n\rw\infty$.\\
We also look at the case $s_k\in\C$. We use the Feller coupling to show that for each $|x|<1$ and $s\in\C$ there exists a random variable $Z_\infty^s(x)$ such that $Z_n^s(x)\xrightarrow{d}Z_\infty^s(x)$ and
$\E{\prod_{k=1}^p Z_n^{s_k}(x_k)}\rw \E{\prod_{k=1}^p Z_\infty^{s_k}(x_k)}$ for $\max_k|x_k|<1$ and $n\rw\infty$.
\end{abstract}

\maketitle
%
%
%
%
%
%
%

\tableofcontents

\newpage
%
%
%
%
\section{Introduction}

The characteristic polynomials of random matrices has been an object of intense interest in recent years.
One of the reasons is the paper of Keating and Snaith \cite{snaith}.
They conjectured that the Riemann zeta function on the critical line could be
modeled by the characteristic polynomial of a random unitary matrix considered on the unit circle.
One of the results in \cite{snaith} is
\begin{theorem}
 Let $x$ be a fixed complex number with $|x|=1$ and $g_n$
be a unitary matrix chosen at random with respect to the Haar measure. Then
$$\frac{\Log\Bigl(\det(I_n-xg_n)\Bigr)}{\sqrt{\frac{1}{2} \log(n)}}
\xrightarrow{d} \mathcal{N}_1 + i\mathcal{N}_2 \text{ for }n\rw\infty$$
and $\mathcal{N}_1,\mathcal{N}_2$ independent, normal distributed random variables.
\end{theorem}
A probabilistic proof of this result can be found in \cite{ashkan}.\\

A surprising fact, proven in \cite{HKOS}, is that a similar result holds for permutation matrices.
A permutation matrix is a unitary matrix of the form $(\delta_{i,\sigma(j)})_{1\leq i,j\leq n}$ with $\sigma\in S_n$
and $S_n$ the symmetric group. It is easy to see that the permutation matrices form a group isomorphic to $S_n$. We call for simplicity both groups $S_n$ and use this identification without mentioning it explicitly. We will, however, use the notation $g\in S_n$ for matrices and $\sigma\in S_n$ for permutations.
We define the characteristic polynomial as
\begin{align}\label{eq_def_Z_n_simple}
Z_n(x)=Z_n(x)(g):=\det(I-x g) \text{ with } x\in\C, g\in S_n.
\end{align}
We can now state the result in \cite{HKOS}:
\begin{theorem}
Let $x$ be a fixed complex number with $|x|=1$, not a root of unity and of finite type. Let $g$
be a $n\times n$ permutation matrix chosen uniformly at random. Then both the real and
imaginary parts of
$$\frac{\Log\Bigl(Z_n(x)\Bigr)}{\sqrt{\frac{\pi}{12}\log(n)}}$$
converge in distribution to a standard normal-random variable.
\end{theorem}
The goal of this paper is to study the moments $\E{Z_n^s(x)}$ of $Z_n(x)$ with respect to the
\emph{uniform measure} on $S_n$ in the three cases
\begin{enumerate}
  \item[(i)]    $|x|<1, s\in\Nr$ (section~\ref{exp Z_n}),
  \item[(ii)]   $|x|<1, s\in\C$  (section~\ref{holoins}),
  \item[(iii)]  $|x|=1, s\in\Nr$ (section~\ref{growth x=1}).
\end{enumerate}
Keating and Snaith have studied in \cite{snaith} also the moments of the characteristic polynomial of a random unitary matrix.
A particular result is
\begin{align}
\mathbb{E}_{U(n)} \left[ \bigl|\det(I-e^{i\theta}U)\bigr|^s \right]
=
\prod_{j=1}^n
\frac{\Gamma(j)\Gamma(j+s)}{(\Gamma(j+s/2))^2}.
\label{eq_moments_Un}
\end{align}
Their proof base on Weyl's integration formula and Selberg's integral.
An alternative proof of \eqref{eq_moments_Un} with representation theory can be found in \cite{bump_gamburd}.
This proof bases on the Cauchy identity (see \cite[chapter 43]{bump}).
%
Unfortunately both arguments does not work for permutation matrices.
We will see that we cannot give a (non trivial) expression for
$\E{Z_n^s(x)}$, but we can calculate the behavior for $n\rw\infty$ and $x,s$ fixed.\\
We use here two types of arguments: a combinatorial one and a probabilistic one. \\
The combinatorial argument (lemma~\ref{lem_cycle_index_theorem}) allows us to write down the generating function of
$\Bigl(\E{Z_n^s(x)}\Bigr)_{n\in\Nr}$ (see theorem~\ref{n_taylor}).
We can use this generating function to calculate the behavior of $\E{Z_n^s(x)}$ for $n\rw\infty$ in the cases (i)
and (iii).
Both calculations rely on the fact that the generating function is a finite product. This is not true anymore for $s\in\C\setminus\Nr$.
We therefore cannot handle case (ii) in the same way.
Here comes in into play the probabilistic argument, namely the Feller coupling (section~\ref{feller2}). The advantage of the Feller coupling is
that it allows us to directly compare $Z_n(x)$ to $Z_{n+1}(x)$.
Of course case (ii) includes case (i), but the calculations in
section~\ref{holoins} are not as beautiful as in section~\ref{exp Z_n}. Another disadvantage of the Feller coupling is that we cannot
use it to calculate the behavior in case (iii).\\
%
The results in section~\ref{exp Z_n} and \ref{holoins} are in fact true for more than one variable.
We introduce in section~\ref{def_def_und_not} a short notation for this and state the theorems in full generality.
We do the proofs only for one or two variables, since they are basically the same as in the general case.
\section{Expectation of $Z_n^s(\bx)$}
\label{exp Z_n}
We give in this section some definitions, write down the generating function for $\E{Z_n^s(\bx)}$
and calculate the behavior of $\E{Z_n^s(\bx)}$ for $n\rw\infty$.
%
%
\subsection{Definitions and notation}
\label{def_def_und_not}
Let us first introduce a short notation for certain functions on defined $\C^p$.
\begin{definition}
\label{def_x_short}
Let $x_1,\cdots,x_p$ be complex numbers.
We write $\bx:=(x_1,\cdots,x_p)$ (similarly for $\bs,\bk,\bn$).
Write $\norm{\bx}:=\max_k(|x_k|)$ for the norm of $\bx$.\\
Let $f:\C^2\rw \C$ and $\bs,\bk\in\C^p$ be given.
Then $f(\bs,\bk):=\prod_{j=1}^p f(s_j,k_j)$.
\end{definition}
The multivariate version of $Z_n(x)$ is
\begin{definition}
We set for $\bx\in \C^p, \bs\in\Nr^p$
\begin{align}Z_n^\bs(\bx)=Z_n^\bs(\bx)(g)
:=
\prod_{k=1}^p\det(I-x_k g)^{s_k}\ \text{ for } g\in S_n\label{defzn}
\end{align}
where $g$ is chosen randomly from $S_n$ with respect to the uniform measure.
\end{definition}

\subsection{An representation for $\E{Z_n^\bs(\bx)}$}
We give now a "good" representation for $Z_n^\bs(\bx)$ and $\E{Z_n^\bs(\bx)}$.
We use here well known definitions and results on the symmetric group.
We therefore restrict ourselves to stating results. More details can be found in \cite[chapter I.2 and I.7]{macdonald}
or in \cite[chapter 39]{bump}.\\

We begin with $Z_n^\bs(\bx)(g)$ for $g$ a cycle.
\begin{example}
\label{bsp_cyc}Let $\sigma=(123\cdots n)$. Then
\begin{align}g=(\delta_{i,\sigma(j)})_{1\leq i,j\leq n}=\left(
                                                             \begin{array}{ccccc}
                                                               0 & 0 & \cdots & 0 & 1 \\
                                                               1 & 0 & \cdots & 0 & 0 \\
                                                               0 & 1 & 0 & \ddots & \vdots \\
                                                               \vdots & \ddots & \ddots & \ddots &  \vdots \\
                                                               0 & \cdots & 0 & 1 & 0 \\
                                                             \end{array}
                                                           \right).\end{align}
Let $\alpha(m):=\exp(m\frac{2\pi i}{n})$. Then
\begin{align}g\left(
                   \begin{array}{c}
                     \alpha(m) \\
                     \alpha(m)^2 \\
                     \vdots \\
                     \alpha(m)^{n-1} \\
                     1 \\
                   \end{array}
                 \right)
                 =
		  \left(
                   \begin{array}{c}
                     1 \\
                     \alpha(m) \\
                         \vdots \\
                     \alpha(m)^{n-2} \\
                     \alpha(m)^{n-1}\\
                   \end{array}
                 \right)
                 =
		\overline{\alpha(m)}
		\left(
                   \begin{array}{c}
                     \alpha(m) \\
                     \alpha(m)^2 \\
                     \vdots \\
                     \alpha(m)^{n-1} \\
                     1 \\
                   \end{array}
                 \right).
\end{align}
Therefore the eigenvalues of $g$ are $\set{\overline{\alpha(m)}} = \set{\alpha(m)}$ and $\det(I-xg)=1-x^n$.\\
\end{example}
%

The arbitrary case now follows easily from example~\ref{bsp_cyc}.
Fix an element $\sigma=\sigma_1 \sigma_2\cdots \sigma_l\in S_n$ with $\sigma_i$ disjoint cycles of length $\la_i$. Then
\begin{align}     g=g(\sigma) =\left(
                                  \begin{array}{cccc}
                                    P_1 &     &      &     \\
                                        & P_2 &      &     \\
                                        &     &\ddots&     \\
                                        &     &      & P_l \\
                                  \end{array}
                                 \right)\end{align}
where $P_i$ is a block matrix of a cycle as in example~\ref{bsp_cyc} (after a possible renumbering of the basis).
Then
\begin{align}
\label{def_zn_la}
Z_n^\bs(\bx)(g)
=
\prod_{k=1}^p\prod_{i=1}^l (1-x_k^{\la_i})^{s_k}.
\end{align}
Next, we look at $\E{Z_n^\bs(\bx)}$. It is easy to see that $Z_n^\bs(\bx)(g)$ only depends on the conjugacy class of $g$.\\

We parameterize the conjugation classes of $S_n$ with partitions of $n$.
\begin{definition}
\label{defpart}
A \emph{partition} $\la$ is a  sequence of nonnegative integers $(\la_1,\la_2,\cdots)$ with  $$\la_1\geq\la_2\geq\cdots \text{ and } \sum_{i=1}^\infty \la_i<\infty.$$
The length $l(\la)$ and the size $|\la|$ of $\la$ are defined as
$$l(\la):=\max\set{i\in\Nr; \la_i\neq0}\text{ and }|\la|:=\sum_{i=1}^\infty \la_i.$$
We set $\la\vdash n:=\set{\la\text{ partition };|\la|=n}$ for $n\in \Nr$.
An element of $\la\vdash n$ is called a \emph{partition of $n$}.
\end{definition}
\textbf{Remark}: we only write the non zero components of a partition.\\
Choose any $\sigma \in S_n$ and write it as $\sigma_1 \sigma_2\cdots \sigma_l$ with $\sigma_i$ disjoint cycles of length $\la_i$. 
Since disjoint cycles commute, we can assume that $\la_1\geq\la_2\cdots\geq \la_l$.
Therefore $\la:=(\la_1,\cdots,\la_l)$ is a partition of $n$.
\begin{definition}
We call the partition $\la$ the \emph{cycle-type} of $\sigma\in S_n$.
\end{definition}
%
%
%
\begin{definition}
Let $\la$ be a partition of $n$. We define $\mathcal{C}_\la\subset S_n$ to be the set of all elements with cycle type $\la$.
\end{definition}
It follows immediately from \eqref{def_zn_la} that $Z_n(\bx)(\sigma)$ only depends on the cycle type of $\sigma$. If
$\sigma,\theta\in S_n$ have different cycle type then $Z_n(\bx)(\sigma)\neq Z_n(\bx)(\theta)$.\\
Therefore two elements in $S_n$ can only be conjugated if they have the same cycle type (since $Z_n(\bx)$ is a class function).
One can find in \cite[chapter 39]{bump} or in \cite[chapter I.7, p.60)]{macdonald} that this condition is sufficient.
The cardinality of each $C_\la$ can be found also in \cite[chapter 39]{bump}.
\begin{lemma}
\label{lem_size_conj}
We have $|\mathcal{C}_\la|=\frac{n!}{z_\la}$ with
\begin{align}z_\la:=\prod_{r=1}^{n} r^{c_r}c_r!
\text{ and } c_r=c_r(\la):=\#\set{i| \la_i=r}.
\end{align}
\end{lemma}
%

We put lemma~\ref{lem_size_conj} and \eqref{def_zn_la} together and get 
\begin{lemma}
 Let $\bx \in\C^p$ and $\bs\in\Nr^p$ be given. Then
\begin{align}
\E{Z_n^{\bs}(\bx)}
=
\sum_{\la \vdash n} \frac{1}{z_\la} \prod_{k=1}^p\prod_{m=1}^{l(\la)}(1-x_k^{\la_m})^{s_k}.
\label{eq_exp_Zn_s_with_la}
\end{align}
\end{lemma}
Obviously, it is very difficult to calculate $\E{Z_n^\bs(\bx)}$ explicitly. It is much easier to write down the generating function.
\subsection{Generating function of $\E{Z^\bs_n(\bx)}$}
\label{sec_gen_fkt}
We now give the definition of a generating function, some lemmas and apply them to
the sequence $\left(\E{Z^\bs_n(\bx)}\right)_{n\in\Nr}$.
\begin{definition}
 Let $(f_n)_{n\in\Nr}$ be given. The formal power series  $f(t):=\sum_{n=0}^\infty f_n t^n$ is called \emph{the generating function} of the sequence $(f_n)_{n\in\Nr}$.
If a formal power series $f(t)=\sum_{n\in\Nr} f_n t^n$ is given then $\left[f\right]_n:=f_n$.
\end{definition}

We will only look at the case $f_n\in\C$ and $f$ convergent.
\begin{lemma}
\label{lem_cycle_index_theorem}
Let $(a_m)_{m\in\Nr}$ be a sequence of complex numbers. Then
\begin{align}
\label{symm fkt}
\sum_{\la} \frac{1}{z_\la} a_\la t^{|\la|}=\exp\left(\sum_{m=1}\frac{1}{m} a_m t^m\right)
\text{ with }
a_\la:=\prod_{i=1}^{l(\la)} a_{\la_i}.
\end{align}

If RHS or LHS of \eqref{symm fkt} is absolutely convergent then so is the other.
\end{lemma}
\begin{proof}
The proof can be found in \cite[chapter I.2, p.16-17]{macdonald} or can be directly verified using the definitions of $z_\la$ and the exponential function.
\end{proof}
%
%
In view of \eqref{eq_exp_Zn_s_with_la} is natural to use lemma~\ref{lem_cycle_index_theorem}
with $a_m = \prod_{j=1}^p (1-x_j^m)^{s_j}$ to write down the generating function of the sequence
$\left(\E{Z^\bs_n(\bx)}\right)_{n\in\Nr}$. We formulate this a lemma, but with $a_m = P(\bx^{m})$ for $P$ a polynomial.
We do this since the calculations are the same as for $a_m = \prod_{j=1}^p (1-x_j^m)^{s_j}$.
\begin{lemma}
\label{generating-allg}
Let $P(\bx)=\sum_{\bk\in\Nr^p} b_{\bk} \bx^\bk$ be a polynomial with $b_{\bk}\in\C$ and $\bx^\bk := \prod_{j=1}^p x_j^{k_j}$.
We set for a partition $\la$
$$
P_\la(\bx):=\prod_{m=1}^{l(\la)}P(\bx^{\la_m})
\text{ with }
\bx^{\la_m} = (x_1^{\la_m}, \cdots, x_p^{\la_m})
\text{ and }
\Omega:=\set{(t,\bx)\subset \C^{p+1};|t|<1, \norm{\bx}\leq 1}.$$
%
%
We then have
\begin{align}\label{generating 1}
\sum_{\la} \frac{1}{z_\la} P_\la(\bx)t^{|\la|}
=
\prod_{\bk\in\Nr^p} (1-\bx^\bk t)^{-b_\bk}
\end{align}
and both sides of \eqref{generating 1} are holomorphic on $\Omega$.\\
We use the principal branch of the logarithm to define $z^s$ for $z\in \C\setminus{\Rr_-}$.
\end{lemma}
%
%
%
%

%
%
\begin{proof}
We only prove the case $p=2$. The other cases are similar.
We use \eqref{symm fkt} with $a_m=P(x_1^m,x_2^m)$ and get
\begin{align*}
\sum_{\la} \frac{1}{z_\la} \prod_{m=1}^{l(\la)}P(x_1^{\la_m},x_2^{\la_m})t^{|\la|}
&=
\exp\left(\sum_{m=1}\frac{1}{m}P(x_1^m,x_2^m) t^m \right)
=
\exp\left(\sum_{k_1,k_2=0}^\infty b_{k_1,k_2} \sum_{m=1}^\infty \frac{t^m}{m} (x_1^{k_1}x_2^{k_2})^m\right)\\
&=
\exp\left(\sum_{k_1,k_2=0}^\infty b_{k_1,k_2} (-1)\Log(1-x_1^{k_1} x_2^{k_2}t)\right)\\
&=
\prod_{{k_1,k_2}=0}^\infty (1-x_1^{k_1}x_2^{k_2}t)^{-b_{k_1,k_2}}.
\end{align*}
The exchange of the sums is allowed since there are only finitely many non zero $b_{k_1,k_2}$.
Since we have used the Taylor-expansion of $\Log(1+z)$ near $0$, we have to assume that $|t x_1^{k_1}x_2^{k_2}|<1$ if $b_{k_1,k_2}\neq 0$.
\end{proof}
We now write down the generating function of $\E{Z_n^{\bs}(\bx)}$.
%
\begin{theorem}
\label{n_taylor}
Let $\bs\in\Nr^p$ and $\bx\in\C^p$. We set
\begin{align}\label{def_fsxt}
f^{(\bs)}(\bx,t)
:=
\prod_{\bk\in\Nr^p}\left(1-\bx^{\bk}t\right)^{-\binom{\bs}{\bk}(-1)^{\bk}}
\text{ with }
\binom{\bs}{\bk}(-1)^{\bk}:= \prod_{j=1}^p \binom{s_j}{k_j}(-1)^{k_j}.
\end{align}
Then
$$
\E{Z_n^{\bs}(\bx)}
=
\left[f^{(\bs)}(\bx,t)\right]_n. $$
\end{theorem}

\begin{proof}
This is \eqref{eq_exp_Zn_s_with_la} and lemma~\ref{generating-allg} with
$P(\bx) = \prod_{j=1}^p (1-x_j)^{s_j}$.\\
\end{proof}
\textbf{Remark}: we use the convention $\E{Z_0^\bs(\bx)}:=1$.\\
\textbf{Remark}: In an earlier draft version, the above proof was based on representation theory. It was at least twice as long and much more
difficult. We wish to acknowledge Paul-Olivier Dehaye, who has suggested the above proof and allowed us to us it.
\begin{corollary}
 For $n\geq 1,p=s=1$, we have that
\begin{align}\label{exp z_n}
\E{Z_n(x)}=1-x.
\end{align}
\end{corollary}
\begin{proof}
$$\frac{1-xt}{1-t}=\sum_{n=0}^\infty t^n-xt\sum_{n=0}^\infty t^n=1+(1-x)\sum_{n=1}^\infty t^n$$
\end{proof}
%
%
%
%
%
%
\subsection{Asymptotics for $\norm{\bx}<1$}
\label{conver x<1}
\begin{theorem}
\label{theo con_s_in_n}
Let $\bs\in\Nr^p$ and $\bx\in\C^p$ with $\norm{\bx}<1$ be fixed. Then
\begin{align}
\lim_{n\rw\infty}\E{Z_n^{\bs}(\bx)}
=\prod_{\bk\in\Nr^p\setminus\set{0}}
\left(1-\bx^{\bk}\right)^{-\binom{\bs}{\bk}(-1)^{\bk}}.
\label{con_s_in_n}
\end{align}
\end{theorem}
%
%
%
There is no problem of convergence in the RHS of \eqref{con_s_in_n} since by definition $\binom{\bs}{\bk}(-1)^{\bk} \neq 0$
only for finitely many $\bk\in\Nr^p\setminus\set{0}$.

There exists several ways to prove this theorem. The proof in this section is a simple, direct computation.
In section~\ref{sec-holo_in_s}, we extend theorem~\ref{theo con_s_in_n} to complex $\bs$ and prove it with probability theory.
The third way to prove this theorem is to use the theorems IV.1 and IV.3 in \cite{FlSe09} which base on function theory.
%
%
We now prove theorem~\ref{theo con_s_in_n} with theorem~\ref{n_taylor} and
\begin{lemma}
Let $a,b\in\Nr$ and $y_1,\cdots,y_a,z_1,\cdots,z_{b}$ be complex numbers with $\max\set{|y_i|,|z_i|} < 1$. Then
$$
\lim_{n\rw\infty} \left[\frac{1}{1-t}\frac{\prod_{i=1}^a(1-y_i t)}{\prod_{i=1}^{b}(1-z_{i} t)}\right]_n
=
\frac{\prod_{i=1}^a(1-y_i)}{\prod_{i=1}^{b}(1-z_{i})}.
$$
\end{lemma}
\begin{proof}
We show this by induction on the number of factors.\\
For $a=b=0$ there is nothing to do, since $\left[\frac{1}{1-t}\right]_n=1$.\\
Induction $(a,b)\rw(a+1,b)$.\\
We set
\begin{align*}
g(t):=\frac{1}{1-t}\frac{\prod_{i=1}^a(1-y_i t)}{\prod_{i=1}^{b}(1-z_{i} t)},
\quad
\gamma:=\frac{\prod_{i=1}^a(1-y_i)}{\prod_{i=1}^{b}(1-z_{i})}.
\end{align*}
%
%
We know by induction that $\lim_{n\rw\infty}[g(t)]_n=\gamma$. We get
%
\begin{align*}
\lim_{n\rw\infty} \left[\frac{1}{1-t}\frac{\prod_{i=1}^{a+1} (1-y_i t)}{\prod_{i=1}^{b}(1-z_{i} t)}\right]_n
&=
\lim_{n\rw\infty}\nth{g(t)(1-y_{a+1} t)}
=
\lim_{n\rw\infty}\nth{g(t)} - \lim_{n\rw\infty} y_{a+1}  \left[{g(t)}\right]_{n-1}\\
&=
(1-y_{a+1})\gamma
=
\frac{\prod_{i=1}^{a+1} (1-y_i)}{\prod_{i=1}^{b}(1-z_{i})}.
\end{align*}
%
%
Induction $(a,b)\rw(a,b+1)$.\\
This case is slightly more difficult. We define $g(t)$ and $\gamma$ as above and write for shortness $z=z_{b+1}$. We have
\begin{align*}
\nth{\frac{1}{1-t}\frac{\prod_{i=1}^{a} (1-y_i t)}{\prod_{i=1}^{b+1}(1-z_{i} t)}}
=
\nth{\frac{g(t)}{(1-zt)}}
=
\nth{g(t)\sum_{k=0}^\infty (zt)^k}
=
\sum_{k=0}^n z^k \left[{g(t)}\right]_{n-k}.
\end{align*}
Let $\eps>0$ be arbitrary. Since $\nth{g(t)}\rw \gamma$, we know that there exists an $n_0\in\Nr$
with $\left|\nth{g(t)}-\gamma \right|<\eps$ for all $n\geq n_0$. We have
\begin{align*}
\sum_{k=0}^n z^k \left[{g(t)}\right]_{n-k}
=
\sum_{k=0}^{n-n_0} z^k \left[{g(t)}\right]_{n-k} +\sum_{k=0}^{n_0-1} z^{n-k}\left[{g(t)}\right]_{k}.
\end{align*}
%
%
Since $|z|<1$, the second sum converges to $0$ as $n\rw\infty$. But
\begin{align*}
\left|\sum_{k=0}^{n-n_0}\gamma z^{k}- \sum_{k=0}^{n-n_0}\left[{g(t)}\right]_{n-k}z^{k} \right|
\leq
|\eps| \sum_{k=0}^{n-n_0}\left| z^{n-k}\right|
\leq
\frac{\eps }{1-|z|}.
\end{align*}
%
%
%
%
%
Since $\eps$ was arbitrary, we are done.
\end{proof}
\begin{corollary}
 \label{con abs in N}
For each $s\in\Nr$ and $x\in\C$ with $|x|<1$, we have
$$\lim_{n\rw\infty} \E{|Z_n(x)|^{2s}}=\prod_{k=1}^s(1-|x|^{2k})^{-\binom{s}{k}^2}
\left|\prod_{0 \leq k_1 <k_2 \leq s} (1-x^{k_1}\overline{x}^{k_2})^{\binom{s}{k_1}\binom{s}{k_2}(-1)^{k_1+k_2+1}}
\right|^2.$$
\end{corollary}
\begin{proof}It follows immediately from \eqref{def_zn_la}, that $\overline{Z_n(x)}=Z_n(\overline{x})$.
We put $p=2$, $s_1=s_2=s$, $x_1=x$ and $x_2=\overline{x}$.
\end{proof}
\textbf{Remark}: the corollary is in fact true for all $s\in\Rr$. One simply has to replace theorem~\ref{theo con_s_in_n} by theorem~\ref{|x|<1} (see section~\ref{holoins}) in the proof.

\section{\label{sec-holo_in_s}Holomorphicity in $s$}
\label{holoins}
%
This section is devoted to extend theorem~\ref{theo con_s_in_n} to $\bs\in\C^p$. We do this by showing that
all functions appearing in theorem~\ref{theo con_s_in_n} are
holomorphic functions in $(\bx,\bs)$ for $\norm{\bx}<1$ and then prove point-wise convergence of these functions.
We do this since we need the holomorphicity in the direct proof of theorem~\ref{|x|<1} (see section~\ref{eq_second_proof_of_main_thm})
and since there are only minor changes between $\bs$ fix and $\bs$ as variables.\\
%
%

We do not introduce here holomorphic functions in more than one variable since we do not need it in
the calculations (except in section~\ref{eq_second_proof_of_main_thm}).
A good introduction to holomorphic functions in more than one variable is
the book ``\emph{From holomorphic functions to complex manifolds}'' \cite{fritzsche}.\\

We now state the main theorem of this section

\begin{theorem}
\label{|x|<1}
We have
\begin{align}
\E{Z_n^{\bs}(\bx)}
\to
\prod_{\bk\in\Nr^p\setminus\set{0}}
\left(1-\bx^{\bk}\right)^{\left(-\binom{\bs}{\bk}(-1)^{\bk}\right)}
\text{ for $n\to \infty$ and all $\bx,\bs \in \C^p$ with $\norm{\bx} < 1$.}
\label{eq_conv_x<1_s_in_C}
\end{align}
%
We use the principal branch of logarithm to define $a^b$ for $a\notin \Rr_{\leq 0}$.
\end{theorem}
%
%
%

%
%
\subsection{Corollaries of theorem~\ref{|x|<1}}
Before we prove theorem~\ref{|x|<1}, we give some corollaries
\begin{corollary}
\label{kor_frac Zn}
We have for $\bs_1,\bs_2,\bx_1,\bx_2\in\C^p$ with $\norm{\bx_1}<1$, $\norm{\bx_2}<1$
$$\E{\frac{Z_n^{\bs_1}(\bx_1)}{Z_n^{\bs_2}(\bx_2)}}
\rw
\prod_{\substack{\bk_1,\bk_2\in\Nr^p\\\bk_1+\bk_2\neq 0}}
\left(1-\bx_1^{\bk_1}\bx_2^{\bk_2}\right)^{-\left(\binom{\bs_1}{\bk_1}\binom{\bs_2+\bk_2-1}{\bk_2}(-1)^{\bk_1}\right)}
\qquad (n\rw\infty).$$
\end{corollary}
\begin{proof}
We use the definition of $\binom{s}{k}$ in \eqref{def_binom} for $s\in\C,k\in\Nr$ (see later).\\
We apply theorem~\ref{|x|<1} for $p':=2p$ and the identity $\binom{-s}{k}=(-1)^k\binom{s+k-1}{k}$.
\end{proof}
\begin{corollary}
 We have for $x_1,x_2,x_3,x_4,s_1,s_2,s_3,s_4\in\C$ with $\max\set{|x_i|}<1$
\begin{align*}
\E{\frac{Z_n^{s_1}(x_1)Z_n^{s_2}(x_2)}{Z_n^{s_3}(x_3)Z_n^{s_4}(x_4)}}
\rw
\prod_{\substack{k_1,k_2,k_3,k_4\in\Nr\\k_1+k_2+k_3+k_4\neq0}}
\left(1-x_1^{k_1}x_2^{k_2}x_3^{k_3}x_4^{k_4}\right)^{\binom{s_1}{k_1}\binom{s_2}{k_2}\binom{s_3+k_3-1}{k_3}\binom{s_4+k_3-1}{k_4}(-1)^{k_1+k_2+1}}
\end{align*}
\end{corollary}
We can also calculate the limit of the Mellin-Fourier-transformation of $Z_n(x)$,
as Keating and Snaith  did in their paper \cite{snaith} for the unitary group.
\begin{corollary}
 We have for $s_1,s_2\in \Rr,x\in\C$ with $|x|<1$
$$
\E{|Z_n(x)|^{s_1}e^{i s_2 \mathrm{arg}(Z_n(x))}}
\to
\prod_{\substack{k_1,k_2\in\Nr\\k_1+k_2\neq0}}^{\infty}
\left(1-x^{k_1} \overline{x}^{k_2}\right)^{\left(\binom{\frac{s_1-s_2}{2}}{k_1}\binom{\frac{s_1+s_2}{2}}{k_2}(-1)^{k_1+1}\right)}.$$
\end{corollary}
\begin{proof}
We have $|z|^{s_1}=z^{s_1/2}\overline{z}^{s_1/2}$ and $e^{i s_2 arg(z)}=\frac{z^{s_2}}{|z|^{s_2}}$.
\end{proof}
%
\subsection{Easy facts and definitions}
We simplify the proof of theorem~\ref{|x|<1} by assuming $p=1$.
We first rewrite $(1-x^m)$ as $r_m e^{i\varphi_m}$ with $r_m>0$ and $\varphi_m\in]-\pi,\pi]$ for all $m\in\Nr$. \\
%
%

\textbf{Convention}: we choose $0<r<1$ fixed and prove theorem~\ref{|x|<1} for $|x|<r$.\\
We restrict $x$ to $\set{|x|<r}$ because some inequalities in the next lemma are only true for $r<1$.
\begin{lemma}
The following hold:
\begin{enumerate}
  \item \label{a_m} $1-r^m\leq r_m\leq 1+r^m$ and $|\varphi_m|\leq \alpha_m$, where $\alpha_m$ is defined in figure~\ref{fig_def_a_m}.
	\begin{figure}[h]
		\centering
 		\includegraphics[width=0.5\textwidth]{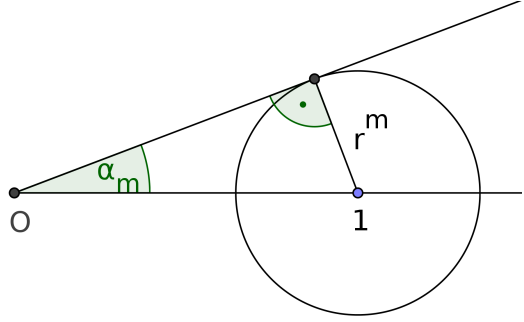}
		\caption{Definition of $\alpha_m$ }
		\label{fig_def_a_m}
	\end{figure}
  \item \label{beta1}One can find a $\beta_1>1$ such that $0\leq\alpha_m\leq \beta_1 r^m$.
  \item \label{beta2}For $-r<y<r$ one can find a $\beta_2=\beta_2(r)>1$ with $|\log(1+y)|\leq \beta_2 |y|$.
  \item \label{beta3}There exists a $\beta_3=\beta_3(r)$, such that for all $m$ and $0\leq y\leq r$
        $$1+y^m\leq \frac{1}{1-y^m}\leq 1+\beta_3 y^m.$$
  \item We have for all $s\in \C$
  $$\Log\left((1-x^m)^{s}\right)\equiv s \Log\left(1-x^m\right)\mod 2\pi i$$
        with $\Log(.)$ the principal branch of logarithm.
\end{enumerate}
\end{lemma}
\begin{proof}
The proof is straight forward. We therefore give only an overview
\begin{enumerate}
 \item We have $|x^m|< r^m$ and thus $1-x^m$ lies inside the circle in figure~\ref{fig_def_a_m}.
 This proves point~\eqref{a_m}
 \item We have that $\sin(\alpha_m) = r^m$ by definition and $\sin(z) \sim z$ for $z \to 0$. This proves point~\eqref{beta1}.
 \item We have $\log(1+y) = \log|1+y|+ i\textrm{arg}(1+y)$.
       This point now follows by takeing a look at $\log|1+y|$ and $\textrm{arg}(1+y)$ separately.
 \item Obvious.
 \item Obvious.

\end{enumerate}
\end{proof}
\begin{lemma}
Let $Y_m$ be a Poisson distributed random variable with $\E{Y_m}=\frac{1}{m}$. Then
         $$\E{y^{(d Y_m)}}=\exp\left(\frac{y^d-1}{m}\right)\text{ for }y,d\geq0.$$
\end{lemma}
%

%
%
\subsection{Extension of the definitions}
\label{ext_def}
The first thing we have to do is to extend the definitions of $Z^s_n(x)$ and $f^{(s)}(x)=\prod_{k=1}^{s}\left(1-x^{k}\right)^{\binom{s}{k}(-1)^{k+1}}$ to holomorphic functions in $(x,s)$.
\subsubsection{Extension of $f^{(s)}(x)$}
We first look at $\binom{s}{k}$. We set
\begin{align}\label{def_binom}
\binom{s}{k}:=\frac{s(s-1)\cdots(s-k+1)}{k!}=\prod_{m=1}^k\frac{s-m+1}{m}.
\end{align}
Obviously, $\binom{s}{k}$ can be extended to a holomorphic function in $s$. \\
%
%
We set
$$
f^{(s)}(x):=\prod_{k=1}^{\infty}\left(1-x^{k}\right)^{\left(\binom{s}{k}(-1)^{k+1}\right)}.
$$
The factor $(1-x^k)^{(..)}$ is well defined, because $\Re(1-x^k)>0$.\\
This definition agrees with the old one, since $\binom{s}{k}=0$ for $k>s$ and $s,k\in\Nr$.
Of course we have to show that the infinite product is convergent.
\begin{lemma}
\label{fsholo}
The function $f^{(s)}(x)$ is a holomorphic function in $(x,s)$.
\end{lemma}
\begin{proof}
Choose $s\in K\subset \C$ with $K$ compact. We have
$$\sup_{s\in K}\left|\frac{s-m+1}{m}\right|\rw 1\text{ for }m\rw \infty.$$
We obtain that for each $a>1$ there exists a $C=C(a,K)$ with $|\binom{s}{k}|\leq C a^k$ for all $k\in\Nr,s\in K$.\\
We choose an $a>1$ such that $ra<1$ and set $\Log(-y):=\log(y)+i\pi$ for $y\in\Rr_{>0}$. Then
\begin{align*}
\sum_{k=1}^\infty\left|\Log\left((1-x^k)^{\left((-1)^{k+1}\binom{s}{k}\right)}\right)\right|
&\leq
\sum_{k=1}^\infty\left|\binom{s}{k}\Log(1-x^k)\right|
=
\sum_{k=1}^\infty\left|\binom{s}{k}\Log(r_k e^{i\varphi_k})\right|\\
&\leq
\sum_{k=1}^\infty\left|\binom{s}{k}\right|(|\Log(1-r^k)|+|\alpha_k|)
\leq
\sum_{k=1}^\infty C(\beta_1+ \beta_2)(ar)^{k}<\infty.
\end{align*}
Since this upper bound is independent of $x$, we can find a $k_0$ such that
$$\mathrm{arg}\left((1-x^k)^{\left((-1)^{k+1}\binom{s}{k}\right)}\right)\in [-\pi/2,\pi/2] \ \forall k\geq k_0\text{ and }|x|<r.$$
Therefore $\sum_{k=k_0}^\infty\Log\left((1-x^k)^{\left((-1)^{k+1}\binom{s}{k}\right)}\right)$ is a holomorphic function.
This proves the holomorphicity of $f^{(s)}(x)$.
\end{proof}
\subsubsection{Extension of $Z^s_n(x)$}
We have found in \eqref{def_zn_la} that
$$
Z_n(x)(g)
=
Z_n^1(x)(g)
=
\prod_{i=1}^{l(\la)} (1-x^{\la_i})
\text{ for }g\in \mathcal{C}_\la.$$

We slightly reformulate this formula.\\
\begin{definition}
\label{def_c_m}Let $\sigma\in S_n$ be given.
We define $C_m=C^{(n)}_m=C^{(n)}_m(\sigma)$ to be the number of cycles of length $m$ of $\sigma$.
\end{definition}
The relationship between partitions and $C_m$ is as follows: if $\sigma\in \mathcal{C}_\la$ is given then
$C_m^{(n)}(\sigma)=\# \set{i; \la_i=m}.$ The function $C^{(n)}_m$ only depends on the cycle-type and is therefore as class function on $S_n$. We get
\begin{align}\label{def_zn_cyc}
Z_n(x)
=
Z_n^1(x)(g)
=
\prod_{m=1}^n(1-x^m)^{C^{(n)}_m}.
\end{align}
Since $\Re(1-x^m)>0$, we can use this equation to extend the definition of $Z_n^s(x)$.
\begin{definition}
\label{z_n_s_complex}
We set for $s\in\C$ and $|x|<1$
$$Z_n^s(x):=\prod_{m=1}^n(1-x^m)^{(sC^{(n)}_m)}.$$
\end{definition}
It is clear that $Z_n^s(x)$ agrees with the old function for $s\in\Nr$
and is a holomorphic function in $(x,s)$ for all values of $C^{(n)}_m$.\\
%
%
%
\subsection{The Feller-coupling}
\label{feller2}
In this subsection we follow the book of Arratia, Barbour and Tavar\'{e} \cite{barbour}. The first thing we mention is
\begin{lemma}
\label{cm_rw_ym}
The random variables $C^{(n)}_m$ converge for each $m\in\Nr$ in distribution to a Poisson-distributed random variable $Y_m$ with $\E{Y_m}=\frac{1}{m}$. In fact, we have for all $b\in\Nr$
$$(C_1^{(n)},C_1^{(n)},\cdots,C^{(n)}_b) \xrightarrow{d} (Y_1,Y_2,\cdots,Y_b) \qquad (n\rw\infty)$$
and the random variables $Y_m$ are all independent.
\end{lemma}
\begin{proof}
See \cite[theorem 1.3]{barbour}.
\end{proof}
Since $Z_n(x)$ and $Z_{n+1}(x)$ are defined on different spaces, it is difficult to compare them. This is where the Feller-coupling comes into play. The Feller-coupling constructs a probability space and new random variables $C_m^{(n)}$ and $Y_m$, which have the same distributions as the $C_m^{(n)}$ and $Y_m$ above and can be compared very well.\\
We give here only a short overview. The details of the Feller-coupling can be found in \cite[p.8]{barbour}.\\

The construction is as follows: let $\left<\xi_t\right>_{t=1}^n$ be a sequence of independent Bernoulli-random variables with $\E{\xi_m}=\frac{1}{m}$. We use the following notation for the sequence
$$
\xi=(\xi_1\xi_2\xi_3\xi_4\xi_5\cdots).
$$
An $m-$spacing is a finite sequence of the form
$$
1\underbrace{0\cdots0}_{m-1\text{ times}} 1
$$
\begin{definition}
Let $C_m^{(n)} = C_m^{(n)}(\xi)$ be the number of m-spacings in $1\xi_2\cdots \xi_n 1$.\\
We define $Y_m = Y_m(\xi)$ to be the number of m-spacings in the whole sequence.
\end{definition}
\begin{theorem}
 \label{feller}We have
          \begin{enumerate}
            \item The above-constructed $C_m^{(n)}(\xi)$ have the same distribution as the $C_m^{(n)}$ in definition~\ref{def_c_m}.
	    \item $Y_m(\xi)$ is a.s. finite and Poisson-distributed with $\E{Y_m}=\frac{1}{m}$.
            \item All $Y_m(\xi)$ are independent.
            \item For fixed $b\in\Nr$ we have
            $$\Pb{\bigl(C_1^{(n)}(\xi) , \cdots , C_b^{(n)}(\xi)\bigr) \neq \bigl(Y_1(\xi) ,\cdots , Y_b(\xi)\bigr)}\rw 0\ (n\rw\infty).$$
          \end{enumerate}
\end{theorem}
\begin{proof}
See \cite[p.8-10]{barbour}.
\end{proof}
We use in the rest of this section only the random variables $C_m^{(n)}(\xi)$ and $Y_m(\xi)$. We
therefore just write $C_m^{(n)}$ and $Y_m$ for them.\\
One might guess that $C_m^{(n)}\leq Y_m$, but this is not true. It is possible that $C_m^{(n)}= Y_m+1$. However, this can only happen if $\xi_{n-m}\cdots\xi_{n+1}=10\cdots 0$. If $n$ is fixed, we have at most one $m$ with $C_m^{(n)}=Y_m+1$. We set
\begin{align}
B_m^{(n)}=\set{\xi_{n-m}\cdots\xi_{n+1}=10\cdots 0}.
\end{align}
\begin{lemma}
\label{weis nicht}
We have
	\begin{enumerate}
          \item $C_m^{(n)}\leq Y_m+\mathbf{1}_{B_m^{(n)}}$,
          \item $\Pb{B_m^{(n)}}=\frac{1}{n+1}$,
          \item $\E{\left|C_m^{(n)}-Y_m\right|}\leq \frac{2}{n+1}$,
          \item $C_m^{(n)}$ does not converges a.s. to $Y_m$.
        \end{enumerate}
\end{lemma}
\begin{proof}
The first point follows immediately from the above considerations.\\
The second point is a simple calculation.\\
The proof of the third point can be found in \cite[p.10]{barbour}. The proof is based on the idea that
$Y_m-C_m^{(n)}$ is (more or less) the number of m-spacings appearing after the $n$'th position.\\
We now look at the last point.
%
%
%
Let $\xi_v\xi_{v+1}\cdots \xi_{v+m+1}=10\cdots01$. We then have for $1\leq{m_0}\leq m-1$ and $v\leq n\leq v+m+1$
$$C_{m_0}^{(n)}=\left\{
              \begin{array}{ll}
                C_{m_0}^{(v)}+1, & \hbox{if } n=v+m_0,\\
                C_{m_0}^{(v)}, & \hbox{if } n\neq v+m_0.
              \end{array}
            \right.
$$
Since all $Y_m<\infty$ a.s. and $\sum_{m=1}^\infty Y_m=\infty$ a.s. we are done.
\end{proof}
%
%
%
%
%
%
%
%
%
\subsection{The limit}
We have $Z_n^s(x)=\prod_{m=1}^n(1-x^m)^{(sC_m^{(n)})}$ and we know $C_m^{(n)} \xrightarrow{d} Y_m$. Does $Z_n^s(x)$ converge in distribution to a random variable? If yes, then a possible limit is
\begin{align}
Z^s_\infty(x):=\prod_{m=1}^\infty\left(1-x^m\right)^{(sY_m)}\label{def z infty}.
\end{align}
We show indeed in lemma~\ref{lem_Zn rw Z_infty} that for $s\in\C$ fixed
\begin{align}\label{eq_Zn rw Z_infty}
Z_n^s(x) \xrightarrow{d} Z_\infty^s(x) \qquad(n\rw\infty).
\end{align}
%
%
%
%
We first show that $Z^s_\infty(x)$ is a good candidate for the limit. We prove that $Z^s_\infty(x)$ and $\E{Z^s_\infty(x)}$ are holomorphic functions in $(x,s)$.\\
Since $Z^s_\infty(x)$ is defined via an infinite product, we have to prove convergence. The following lemma supplies us with the needed tools.
\begin{lemma}
\label{vorbez<1}We have:
	\begin{enumerate}
          \item $\sum_{m=1}^\infty Y_m r^m$ is a.s. absolute convergent for $|r|<1$.
          \item $\E{\prod_{m=1}^\infty(1\pm r^m)^{(\sigma Y_m)}}$ is finite for all $\sigma\in \Rr$.
          \item $\E{\prod_{m=1}^\infty\exp(t\alpha_mY_m)}$ is finite for all $t\in \Rr$.
        \end{enumerate}
\end{lemma}
\begin{proof} \ud{First part:}\\
we prove the absolute convergence of the sum $\sum_{m=1}^\infty Y_m r^m$ by showing
$$
\limsup_{m\rw\infty}\sqrt[m]{|Y_m r^m|}<1 \ \text{a.s.}
$$
We fix an $a$ with $r<a<1$ and set $A_m:=\set{\ Y_mr^m>a^m}=\set{\sqrt[m]{|Y_m r^m|}>a}$. Then
\begin{align*}
\Pb{\limsup_{m\rw\infty}\sqrt[m]{|Y_m r^m|}<1}
&\geq
1-\Pb{\limsup_{m\rw\infty}\sqrt[m]{|Y_m r^m|}>a}
=
1-\Pb{\cap_{n=1}^\infty \cup_{m=n}^\infty A_m}\\
&=
1-\Pb{\limsup (A_m)}.
\end{align*}
We get with Markov's inequality
$$
\Pb{Y_m>\left(\frac{a}{r}\right)^m}
\leq
\frac{\E{Y_m}}{\left(\frac{a}{r}\right)^m}=\frac{1}{m}\left(\frac{r}{a}\right)^m.
$$
Therefore $\sum_m \Pb{Y_m>\left(\frac{a}{r}\right)^m}<\infty$. It follows from the Borel-Cantelli-lemma that
$$\Pb{\limsup (A_m)}=0.$$
\ud{Second part:}\\
\ud{Case $\sigma>0$}.\\
We only have to look at the terms with a plus, since the $(1-r^m)^{(\sigma Y_m)}\leq 1$.\\
We get with monotone convergence
\begin{align*}
\E{\prod_{m=1}^\infty(1+ r^m)^{(\sigma Y_m)}}
&=
\lim_{m_0\rw\infty}\E{\prod_{m=1}^{m_0}(1+r^m)^{(\sigma Y_m)}}\\
&= \lim_{m_0\rw\infty}\prod_{m=1}^{m_0}\exp\left(\frac{(1+r^m)^\sigma-1}{m}\right)
\leq \lim_{m_0\rw\infty}\prod_{m=1}^{m_0}\exp\left(\frac{(1+r^m)^{\lceil\sigma\rceil}-1}{m}\right).
\end{align*}
We have for $\lceil\sigma\rceil$ fixed and $m$ big enough $$\bigl((1+r^m)^{\lceil\sigma\rceil}-1\bigr)=r^m\left(\lceil\sigma\rceil+\sum_{k=2}^{\lceil\sigma\rceil}\binom{\lceil\sigma\rceil}{k}(r^m)^{k-1}\right)\leq 2\lceil\sigma\rceil r^m$$
It follows
\begin{align*}
\lim_{m_0\rw\infty}\prod_{m=1}^{m_0}\exp\left(\frac{(1+r^m)^\sigma-1}{m}\right)\leq C\prod_{m=1}^{\infty}\exp(2\lceil\sigma\rceil r^m)<\infty.
\end{align*}
The constant $C$ only depends on $r$ and $\lceil\sigma\rceil$.\\
\ud{Case $\sigma<0$} \\
We only have to look at the terms with a minus sign. We have $$(1-r^m)^{(\sigma Y_m)}=\frac{1}{(1-r^m)^{(|\sigma| Y_m)}}\leq (1+\beta_3 r^m)^{(|\sigma| Y_m)}$$
We can now argue as in the case $\sigma>0$.\\
\ud{Third part:}\\
the product $\prod_{m=1}^\infty\exp(t\alpha_mY_m)$ is a.s. well defined, since $$\Log\left(\prod_{m=1}^\infty\exp(t\alpha_mY_m)\right)
=\sum_{m=1}^\infty t\alpha_mY_m
\leq t\sum_{m=1}^\infty Y_m \beta_1 r^m
<\infty \text{ a.s.} .$$
We get
\begin{align*}
\E{\prod_{m=1}^\infty\exp(t\alpha_mY_m)}
&\leq
\E{\prod_{m=1}^\infty\exp(t\beta_1 r^m Y_m)}
=
\prod_{m=1}^\infty\exp\left(\frac{e^{t\beta_1 r^m}-1}{m}\right).
\end{align*}
Since $t\beta_1 r^m$ is bounded, we can find a constant $C=C(r,t)$ with $e^{t\beta_1 r^m}-1\leq C t\beta_1 r^m$. We get
\begin{align*}
\E{\prod_{m=1}^\infty\exp(t\alpha_mY_m)}
&\leq
 \prod_{m=1}^\infty\exp\left(\frac{C(t\beta_1 r^m)}{m}\right)
\leq
\exp\left(\sum_{m=1}^\infty Ct\beta_1 r^m\right)<\infty.
\end{align*}
\end{proof}
Now, we can prove:
\begin{lemma}
\label{zs_holo}
$Z_\infty^s(x)$ is a.s. a holomorphic function in $(x,s)$.
\end{lemma}
\begin{proof} We have
\begin{align*}
\sum_{m=1}^\infty \left|\Log\left((1-x^m)^{(sY_m)}\right)\right|
&\leq
\sum_{m=1}^\infty |sY_m||\Log(1-x^m)|
\leq
|s|\sum_{m=1}^\infty Y_m(|\log(r_m)|+|i\varphi_m|)\\
& \leq
|s|\sum_{m=1}^\infty Y_m(\beta_1+\beta_2)r^m
<\infty\hbox{ a.s. by Lemma \ref{vorbez<1}}.
\end{align*}
Since the sum is a.s. finite, we can find a.s. an $m_0$, such that $\arg \left((1-x^m)^{(sY_m)}\right)\in [-\pi/2,\pi/2]$ for all $m\geq m_0$. Then $\sum_{m=m_0}^\infty\Log\left((1-x^m)^{(sY_m)}\right)$ is a holomorphic function and so is $Z^s_\infty(x)$.
\end{proof}
%
%
%
%
%
%
%
%
%
%
%
%
%
\begin{lemma}
\label{z<1_holo}
All moments of $Z_\infty^s(x)$ exist. $\E{Z_\infty^s(x)}$ is holomorphic in $(x,s)$ with $\E{Z_\infty^s(x)}=f^{(s)}(x).$
\end{lemma}
\begin{proof}
Let $s\in K\subset \C$ with $K$ compact.\\
\ud{Step 1}: existence of $\E{Z_\infty^s(x)}$.\\
Let $s=\sigma+it$. Then
\begin{align*}
|Z_\infty^s(x)|
&=
\prod_{m=1}^\infty \left|(1-x^m)^{(sY_m)}\right|
=
\prod_{m=1}^\infty r_m^{(\sigma Y_m)}e^{-tY_m\varphi_m}
\leq
\left(\prod_{m=1}^\infty (1+\beta_3 r^m)^{(|\sigma| Y_m)}\right) \left(\prod_{m=1}^\infty e^{tY_m \alpha_m}\right).
\end{align*}
We set $\sigma_0:=\sup_{s\in K}|\sigma|$ and $t_0:=\sup_{s\in K}|t|$.
We define
\begin{align}
F(r)=F(r,K)
:=\left(\prod_{m=1}^\infty (1+\beta_3 r^m)^{(\sigma_0 Y_m)}\right) \left(\prod_{m=1}^\infty e^{t_0Y_m\alpha_m}\right).\label{F(r)}
\end{align}
It follows from the Cauchy Schwarz inequality (for $L^2$) and lemma~\ref{vorbez<1}, that $\E{F(r)}$ is finite. Therefore $\E{Z_\infty^s(x)}$ exists.\\
%
%
%
\ud{Step 2}: the value and the holomorphicity of $\E{Z_\infty^s(x)}$.\\
We have
\begin{align}
\E{Z_\infty^s(x)}
&=
\E{\prod_{m=1}^{\infty}\left(1-x^m\right)^{sY_m}}
=
\prod_{m=1}^{\infty}\E{(1-x^m)^{sY_m}}
=
\prod_{m=1}^{\infty}\exp\left(\frac{(1-x^m)^s-1}{m}\right)\nonumber\\
&=
\exp\left(\sum_{m=1}^\infty\frac{(1-x^m)^s-1}{m}\right).
\end{align}
The exchange of the product and $\E{..}$ is justified by step 1. The exchange of $\exp$ and the product is justified by the following calculation.\\
We need Newton's binomial series to calculate the last expression. We have
\begin{align}
(1+x)^s=\sum_{k=0}^\infty\binom{s}{k}x^k\text{ for }|x|<1,s\in\C
\end{align}
where the sum on the right side is absolutely convergent (See \cite[p.26]{busam}). \\
We get
\begin{align*}
\sum_{m=1}^\infty\frac{(1-x^m)^s-1}{m}
&=
\sum_{m=1}^\infty\frac{1}{m}\sum_{k=1}^\infty\binom{s}{k}(-x^m)^k
=
\sum_{k=1}^\infty(-1)^k\binom{s}{k}\sum_{m=1}^\infty\frac{(x^k)^m}{m}\\
&=
\sum_{k=1}^\infty(-1)^{k+1}\binom{s}{k}\Log(1-x^k).
\end{align*}
We have to justify the exchange of the two sums in the first line.\\
We have seen in the proof of lemma~\ref{fsholo} that for each $a>1$ there exists a constant $C=C(a,K)$ with $\left|\binom{s}{k}\right| <C a^k$ for $k\in \Nr,s\in K$. Now choose $a>1$ with $ar<1$. We get

\begin{align*}
\sum_{m=1}^\infty\sum_{k=1}^\infty \left|\frac{1}{m}\binom{s}{k}(-x^m)^k\right|
&\leq
\sum_{m=1}^\infty\sum_{k=1}^\infty \frac{1}{m} Ca^k(r^m)^k
=
C\sum_{m=1}^\infty\frac{1}{m}\frac{ar^m}{1-ar^m}<\infty.
\end{align*}
\ud{Step 3}: holomorphicity of $\E{Z_\infty^s(x)}$.\\
We know from step 2 and the definition of $f^{(s)}(x)$ that
$$\E{Z_\infty^s(x)}
=
f^{(s)}(x).$$
Since we have shown in lemma~\ref{fsholo} the holomorphicity of $f^{(s)}(x)$, we are done.\\
\ud{Step 4}: existence of the moments.\\
We have $\left(Z_n^s(x)\right)^m=Z_n^{(ms)}(x)$ for each $m\in\Nr$ and $\overline{Z_n^s(x)}=Z_n^{\overline{s}} (\overline{x})$ (use $\overline{e^z}=e^{\overline{z}}$).\\
The existence of the moments now follows from step 1.
\end{proof}

%
%
%
%
%
%
%
%
%
\subsection{Convergence to the limit}
We have so far extended the definitions and found a possible limit. We complete the proof of theorem~\ref{|x|<1} by showing that
$$\E{Z^s_n(x)}\rw \E{Z^s_\infty(x)}=f^{(s)}(x)\qquad (n\rw\infty)$$
for $|x|<r$ and $s\in\C$.\\
We give here two different proofs. The idea of the first is to prove $Z_n(x)\xrightarrow{d}Z_\infty(x)$ and
then use uniform integrability (see lemma~\ref{lem_uniform integrabel}). The idea of the second is to prove theorem~\ref{|x|<1} for $x,s\in\Rr$ and then apply the theorem of Montel.\\
Note that the second proof does not imply $Z_n(x)\xrightarrow{d}Z_\infty(x)$. One would need that $Z_\infty(x)$ is uniquely defined by its moments. Unfortunately we have not been able to prove or disprove this.
%
%
%

\subsubsection{First proof of theorem~\ref{|x|<1}}
\begin{lemma}
\label{lem_Zn rw Z_infty}We have for all fixed $x$ and $s$
\begin{align}
s\sum_{m=1}^n C_m^{(n)}\Log(1-x^m)
&\xrightarrow{d}
 s\sum_{m=1}^\infty Y_m \Log(1-x^m) \qquad (n\rw\infty),\\
Z_n^{s}(x)
&\xrightarrow{d}
 Z_\infty^{s}(x) \qquad (n\rw\infty).
\end{align}
\end{lemma}
\begin{proof}
Since the exponential map is continuous, the second part follows immediately from the first part.\\
We know from lemma~\ref{weis nicht} that
\begin{align}
\E{\left|C_m^{(n)}-Y_m\right|}\leq \frac{2}{n+1}.
\end{align}
We get
\begin{align*}
\E{\left|s\sum_{m=1}^n (Y_m-C_m^{(n)})\Log(1-x^m)\right|}
&\leq
|s|\sum_{m=1}^n\E{|Y_m-C_m^{(n)}|}|\Log(1-x^m)|\\
&\leq
\frac{2|s|}{n+1}\sum_{m=1}^n (\beta_1+\beta_2)r^m
=
\frac{2|s|r}{n+1}\frac{\beta_1+\beta_2}{1-r}\longrightarrow 0 \quad (n\rw\infty).
\end{align*}
\end{proof}

Weak convergence does not imply automatically convergence of the moments. One needs some additional properties.
We introduce therefore
\begin{definition}
 A sequence of (complex valued) random variables $(X_m)_{m\in\Nr}$ is called uniformly integrable if
$$\sup_{n\in\Nr} \E{|X_n|\mathbf{1}_{|X_n|>c}}\longrightarrow 0 \text{ for } c\to \infty.$$
\end{definition}
\begin{lemma}
\label{lem_uniform integrabel}
Let $(X_m)_{m\in\Nr}$ be uniformly integrable and assume that $X_n\xrightarrow{d} X$. Then
$$\E{X_n} \longrightarrow \E{X}.$$
\end{lemma}
\begin{proof}
See \cite[chapter 6.10, p. 309]{gut}.
\end{proof}
We  finish the proof of theorem~\ref{|x|<1}. Let $s\in K \subset \C$ with $K$ compact.
We have found in the proof of lemma~\ref{z<1_holo}
\begin{align}
|Z_\infty^s(x)|
&\leq
F(r)
=
\left(\prod_{m=1}^\infty (1+\beta_3 r^m)^{(|\sigma_0| Y_m)}\right) \left(\prod_{m=1}^\infty e^{t_0Y_m \alpha_m}\right)
\label{eq_upper_bound_for_Zn}
\end{align}
with $\E{F(r)}<\infty$. It is possible that $C_m^{(n)}=Y_m+1$ and so the inequality $|Z_n^s(x)|<F(r)$ does not have to be true. We replace therefore $Y_m$ by $Y_m+1$ in the definition of $F(r)$. This new $F(r)$ fulfills $\E{F(r)}<\infty$ and
$$|Z_n^s(x)|<F(r) \ , |Z_\infty^s(x)|<F(r) \ \forall n\in\Nr.$$
We get
$$
\sup_{n\in\Nr} \E{|Z^s_n(x)|\mathbf{1}_{|Z^s_n(x)|>c}}
\leq
\E{F(r)\mathbf{1}_{F(r)>c}} \longrightarrow 0 \qquad (c\rw\infty).
$$
The sequence $(Z_n^s(x))_{n\in\Nr}$ is therefore uniformly integrable and theorem \ref{|x|<1} follows immediately from lemma~\ref{lem_Zn rw Z_infty} and \ref{lem_uniform integrabel}.

\subsubsection{Second proof of theorem \ref{|x|<1}}
\label{eq_second_proof_of_main_thm}

We first prove the convergence in the special case $x\in [0,r[, s\in [1,2]$.
\begin{lemma}
\label{convergence1}
For all $x\in [0,r[$ and $s\in [1,2]$, we have
$$\E{Z^s_n(x)}\rw \E{Z^s_\infty(x)}\qquad (n\rw\infty).$$
\end{lemma}
%
\begin{proof}Let $x$ and $s$ be fixed.\\
''$\leq$'' Let $m_0$ be arbitrary and fixed. For $n\geq m_0$ we have
$$\E{Z^s_n(x)}
=\E{\prod_{m=1}^n\left(1-x^m\right)^{sC_m^{(n)}}}
\leq \E{\prod_{m=1}^{m_0}\left(1-x^m\right)^{sC_m^{(n)}}}$$

We know from lemma~\ref{cm_rw_ym} that
$$(C_1^{(n)},C_1^{(n)},\cdots,C^{(n)}_{m_0}) \xrightarrow{d} (Y_1,Y_2,\cdots,Y_{m_0}) \qquad (n\rw\infty).$$
The function $\prod_{m=1}^{m_0}(1-x^m)^{(s \cdot c_m)}$ is clearly continuous in $(c_1,\cdots, c_{m_0})$ and bounded by $1$. We therefore get
$$\E{\prod_{m=1}^{m_0}\left(1-x^m\right)^{sC_m^{(n)}}}\rw \E{\prod_{m=1}^{m_0}\left(1-x^m\right)^{sY_m}}\qquad (n\rw\infty).$$
Since $m_0$ was arbitrary and $(1-x^m)^{sY_m}\leq 1$, it follows with dominated convergence that
$$
\limsup_{n\rw\infty}\E{Z^s_n(x)}
\leq
\inf_{m_0}\E{\prod_{m=1}^{m_0}\left(1-x^m\right)^{sY_m}}
=
\E{\prod_{m=1}^{\infty}\left(1-x^m\right)^{sY_m}}=\E{Z^s_\infty(x)}.
$$
''$\geq$'' The second part is more difficult. Here we need the Feller coupling.\\
Remember: $B_m^{(n)}=\set{\xi_{n-m}\cdots\xi_{n+1}=10\cdots 0}$.\\
If $C_m^{(n)}\leq Y_m$ a.s. we would not have any problems. But it may happen that $C_m^{(n)}= Y_m+1$.
If $\xi_{n+1}=1$, then $C_m^{(n)}\leq Y_m$. For this reason we have
\begin{align} \label{zn cond n-1}
Z_n(x)
=\sum_{m=0}^n Z_n(x)\mathbf{1}_{B^{(n)}_m}
=Z_{n}(x)\mathbf{1}_{B^{(n)}_0} + \sum_{m=1}^n (1-x^m)Z_{n-m}(x)\mathbf{1}_{B^{(n)}_m}.
\end{align}
We choose $1>\eps>0$ fixed and $m_0$ large enough, such that $(1-x^m)^s>1-\eps$ for $m\geq m_0$. Then
\begin{align*}
\E{Z^s_n (x)}
&\geq \E{\sum_{m=1}^n (1-x^m)^sZ^s_{n-m}(x)\mathbf{1}_{B_m}}
\geq \E{\sum_{m=1}^n (1-x^m)^sZ^s_{\infty}(x)\mathbf{1}_{B_m}}\\
&\geq \E{\sum_{m=m_0+1}^\infty(1-\eps)Z_\infty(x)\mathbf{1}_{B_m}}+\E{\sum_{m=1}^{m_0}(1-x)Z_\infty^s(x)\mathbf{1}_{B_m}}.
\end{align*}
The last summand goes to $0$ by the Cauchy Schwarz inequality, since $\Pb{B^{(n)}_m}=\frac{1}{n+1}$ and $Z_\infty^s(x)$ has finite expectation for all $s$. We can replace $m_0$ by $0$ in the other Summand with the same argument.
\end{proof}
%
%
%
We have proven theorem \ref{|x|<1} for some special values of $x$ and $s$. This proof is based on the fact that $0<(1-x^m)\leq 1$ for $x\in [0,1[$.
Therefore we cannot use this proof directly for arbitrary $(x,s)$.
We could try to modify the proof, but this turns out to be rather complicated.
An easier way is to use the theorem of Montel (See \cite[p.230]{busam} or \cite[p.23]{fritzsche}).\\
Suppose that there exists a $(x_0,s_0)$ and a subsequence $\Lambda$ such that\\
$|\E{Z_n^{s_0}(x_0)}-\E{Z_\infty^{s_0}(x_0)}|>\eps$ for some $\eps$ and all $n\in\Lambda$.
We have found in the first proof in \eqref{eq_upper_bound_for_Zn} (and in the proof of lemma~\ref{z<1_holo}) an upper bound $F(r)$
for the sequence $\left(\E{Z_n^{s}(x)}\right)_{n\in\Lambda}$ for $|x|<r$.
The sequence $\left(\E{Z_n^{s}(x)}\right)_{n\in\Lambda}$ is thus locally bounded and we can use the theorem of Montel.
We therefore can find a subsequence $\Lambda'$ of $\Lambda$ and a holomorphic function
$g$ with $\E{Z_n^s(x)}\rw g$ on $B_r(0)\times K$ (for $n\in\Lambda'$). But $g(x)$ has to agree with $\E{Z_\infty^s(x)}$ on $[0,r[\times[1,2]$. Therefore $g(x)=\E{Z_\infty^s(x)}$. But this is impossible since $g(x_0,s_0)\neq \E{Z_\infty^{s_0}(x_0)}$.
This completes the proof of theorem~\ref{|x|<1}.\qed

\section{Growth Rates for $|x|=1$}
\label{growth x=1}
We consider in this section only the case $p=2$ and $x=x_1=\overline{x_2}$.
We assume $s_1,s_2\in\Nr$, $|x|=1$ and $x$ not a root of unity, i.e $x^k\neq1$ for all $k\in\Zr\setminus\set{0}$.\\
We first calculate the growth rate of $\E{Z_n^{s_1}(x)Z_n^{s_2}(\overline{x})}$ for $s_2=0$ (see lemma \ref{growth s2=0}
and \ref{lemma s=4m+2}) and then for $s_2$ arbitrary (see theorem \ref{thm_z^s1 z^s2 groth}).\\
The main results in this section can be obtained by using the theorems VI.3 and VI.5 in \cite{FlSe09}.
This theorems base on Cauchy's integral formula and are very general.
To apply them, one just has to take a look at the singularities of the generating function.
We do not use this theorems since we can reach our target with
a much simpler tool: a partial fraction decomposition.
The advantage is that we get the asymptotic behavior with a simple calculation
(see \eqref{partfracdecompo}, \eqref{geomtr_ord_l} and \eqref{growth_rate}) and
the calculations after \eqref{growth_rate} are unchanged.
%
\subsection{Case $s_2=0$}
We use the generating function of theorem \ref{n_taylor}. We have to calculate the growth rate of $\left[f^{(s)}(x,t)\right]_n$.
But $f^{(s)}(x,t)$ is a quite complicated expression and we therefore express it in a different way.
Since $f^{(s)}(x,t)$ is a rational function, we can do a partial fractional decomposition with respect to $t$ (and $x$ fixed).
\begin{align}
\label{partfracdecompo}
f^{(s)}(x,t)
=
\prod_{k=0}^s\left(1-x^kt\right)^{(-1)^{k+1}\binom{s}{k}}
=
P(t)+\sum_{k\text{ even}} \sum_{l=1}^{\binom{s}{k}} \frac{a_{k,l}}{(1-x^kt)^l}
\end{align}
 where $P$ is a polynomial and $a_{k,l}\in \C$.
This formulation can be found in \cite[chapter 69.4, p.401]{Heuser1993}.
Note that at this point we need the condition that $x$ is not a root of unity. If $x$ is a root of unity, some
factors can be equal and cancel or increase the power. For example, we have $f^{(s)}(1,t)=1$ for all $s\in\C$.\\

What is the growth rate of $\left[\frac{1}{(1-x^kt)^l}\right]_n$? We have for $l\in\Nr$ and $|t|<1$
\begin{align}
\label{geomtr_ord_l}
\frac{1}{(1-t)^l}
=
\frac{1}{(l-1)!}\sum_{n=0}^\infty \Bigl(\prod_{k=1}^{l-1}(n+k)\Bigr)t^n.
\end{align}
This equation can be shown by differentiating the geometric series.
We get
$$
\left[ \frac{1}{(1-x^k t)^l}\right]_n \sim \frac{n^{l-1}(x^k)^n}{(l-1)!}.
$$
Recall that $A(n)\sim B(n)$ if $\lim_{n\rw\infty}\frac{A(n)}{B(n)}=1$. Since $|x|=1$, we have
\begin{align}
\label{growth_rate}
\left|\left[\frac{1}{(1-x^kt)^l}\right]_n\right|
\sim
\frac{n^{l-1}}{(l-1)!}.
\end{align}
We know from \eqref{growth_rate} the growth rate of each summand in \eqref{partfracdecompo}.
Since we have only finitely many summands, only these $\frac{a_{k,l}}{(1-x^kt)^l}$ are relevant with $l$ maximal and $a_{k,l}\neq 0$.
The LHS of \eqref{partfracdecompo} has for $t = \overline{x}^k$ a pole of order $\binom{s}{k}$ (for $0 \leq k \leq s$ and $k$ even).
We therefore have $a_{k,\binom{s}{k}} \neq 0$ since there is only one summand on the RHS of \eqref{partfracdecompo}
with a pole at $t = \overline{x}^k$ of order at least $\binom{s}{k}$.
Before we can write down the growth rate of $\E{Z_n^s(x)}$, we have to define
\begin{definition}
\label{def_C1}
Let $s,k_0\in\Nr$ with $0\leq k_0\leq s$ and $k_0$ even. We set
\begin{align}
\label{eq_def C(k0)}
C(k_0)
:=
\frac{1}{\bigl(\binom{s}{k_0}-1\bigr)!}
\prod_{k\neq k_0}\left(1-x^k \overline{x}^{k_0}\right)^{(-1)^{k+1}\binom{s}{k}}.
\end{align}
\end{definition}
We put everything together and get
\begin{lemma}
 \label{growth s2=0}
We have for $s\neq 4m+2$
\begin{align}\label{eq E(z_n^s) growth rate}
\E{Z_n^{s}(x)}
\sim
  n^{\binom{s}{k_0}-1} C(k_0)(x^{k_0})^n
\end{align}
with 
$$
k_0:=
\left\{
  \begin{array}{ll}
    2m,  & \hbox{for }s=4m, \\
    2m,  & \hbox{for }s=4m+1, \\
    2m+2,& \hbox{for }s=4m+3.
  \end{array}
\right.
$$
%
%
%
%
\end{lemma}
\begin{proof}
We have to calculate $M:=\max_{k\text{ even}}\binom{s}{k}$.
A straight forward verification shows that $M=\binom{s}{k_0}$ and that there is only one summand with exponent $M$ in the case $s\neq 4m+2$. We apply \eqref{growth_rate} and get
$$\E{Z_n^{s}(x)}
\sim
  n^{\binom{s}{k_0}-1} \frac{a_{k_0,\binom{s}{k_0}}}{\bigl(\binom{s}{k_0}-1\bigr)!}(x^{k_0})^n.$$
By the residue theorem
$a_{k_0,\binom{s}{k_0}}
=
\prod_{k\neq k_0}\left(1-x^k \overline{x}^{k_0}\right)^{(-1)^{k+1}\binom{s}{k}}.$
This proves \eqref{eq E(z_n^s) growth rate}.\\
%
%
%
%
\end{proof}
We see in the next lemma that there can appear more than one constant.
This is the reason why we write $C(k_0)$ for the constant and not $C$ or $C(s,x)$.\\
The case $s=4m+2$ is a little bit more difficult, since there are two maximal terms, i.e. $\binom{4m+2}{2m}=\binom{4m+2}{2m+2}$.
%
\begin{lemma}
\label{lemma s=4m+2}If $s=4m+2$ then
\begin{align}\label{s=4m+2}
\E{Z_n^{s}(x)}
\sim
n^{\binom{4m+2}{2m}-1} \Bigl(C(2m)(x^{2m})^n+C(2m+2) (x^{2m+2})^n\Bigr)
\end{align}
with $C(2m)(x^{2m})^n+C(2m+2)(x^{2m+2})^n=0$ for at most one $n$.\\
\end{lemma}
\begin{proof}
A straight forward verification as in lemma~\ref{growth s2=0} shows that $M=\binom{4m+2}{2m}=\binom{4m+2}{2m+2}$.
Now we have two summands with a maximal $l$.\\
To prove \eqref{s=4m+2}, we have to show $C(2m)(x^{2m})^n+C(2m+2)(x^{2m+2})^n=0$ for only finitely many $n$.
But $C(2m)(x^{2m})^n+C(2m+2)(x^{2m+2})^n=0$ implies $x^{2n}=-\frac{C(2m)}{C(2m+2)}$.
Since $x$ is not a root of unity, all $x^k$ are different.
%
\end{proof}
%
%
%
%
%
%
%
%
%
%
%
%
%
\subsection{Case with $s_2$ arbitrary}
We argue as before.\\
Some factors appearing in $f^{(s,s)}(x,\overline{x},t)$ (see \eqref{def_fsxt}) are equal, so we have to collect them before we can write down the partial fraction decomposition.\\
\begin{align}\label{partfracdecompo2}
f^{(s_1,s_2)}(x,\overline{x},t)=\prod_{k_1=0}^{s_1}\prod_{k_2=0}^{s_2} \left(1-x^{k_1}\overline{x}^{k_2}t\right)^{(-1)^{k_1+k_2+1}\binom{s_1}{k_1}\binom{s_2}{k_2}}
=\prod_{k=-s_2}^{s_1}(1-x^kt)^{S(k)}
\end{align}
with $$S(k)=\sum_{j=0}^\infty \binom{s_1}{k+j}\binom{s_2}{j}(-1)^{k+2j+1}.$$
To calculate $S(k)$ explicit, we need Vandermonde's identity for binomial coefficients:
\begin{align}
\binom{m_1+m_2}{q}=\sum_{j=0}^{\infty}\binom{m_1}{q-j}\binom{m_2}{j}.
\end{align}
We get with $m_1=s_1,m_2=s_2, q=s_1-k$ and $\binom{m}{k}=\binom{m}{m-k}$ that
\begin{align}
S(k)=(-1)^{k+1}\binom{s_1+s_2}{s_1-k}
\end{align}
and therefore
$$
f^{(s_1,s_2)}(x,\overline{x},t)
=
\prod_{k=-s_2}^{s_1}(1-x^kt)^{(-1)^{k+1}\binom{s_1+s_2}{s_1-k}}
=
\prod_{k=-s_2}^{s_1}(1-x^kt)^{(-1)^{k+1}\binom{s_1+s_2}{s_2+k}}.
$$
Before we look at the growth rate of $\E{Z_n^{s_1}(x)Z_n^{s_2}(\overline{x})}$, we define
\begin{definition}
\label{def_C2}
We set for $s_1,s_2\in\Nr,k_0\in \Zr$ with $-s_2\leq k_0 \leq s_1$ and $k_0$ even
\begin{align}\label{eq_def_C(k_0)}
C(k_0)
=
C(s_1,s_2,k_0,x)
=
\frac{1}{\Bigl(\binom{s_1+s_2}{s_2+k_0}-1\Bigr)!}
\prod_{k\neq k_0}(1-x^k \overline{x}^{k_0})^{(-1)^{k+1}\binom{s_1+s_2}{s_2+k}}.
\end{align}
\end{definition}
\textbf{Remark}: definition~\ref{def_C1} is a special case of definition \ref{def_C2}.
Therefore there is no danger of confusion and we can write $C(k_0)$ for both of them.\\

We get
\begin{theorem}
\label{thm_z^s1 z^s2 groth}
\hfill
\begin{itemize}
  \item If $s_1-s_2\neq 4m+2$ then
    \begin{align}\label{eq E(z_n^s1 z_n^s2) growth rate1}
       \E{Z_n^{s_1}(x)Z_n^{s_2}(\overline{x})}
         \sim
       n^{\binom{s_1+s_2}{\lfloor (s_1+s_2)/2\rfloor}-1} C(k_0) (x^{k_0})^n
    \end{align}
    with $k_0:=\left\{
      \begin{array}{ll}
        \frac{s_1-s_2}{2} , & \hbox{for } s_1-s_2=4m,\\
        \frac{s_1-s_2-1}{2}, & \hbox{for } s_1-s_2=4m+1,\\
        \frac{s_1-s_2+1}{2}, & \hbox{for } s_1-s_2=4m+3.
      \end{array}
      \right.$
  \item If $s_1-s_2 = 4m+2$ we set $k_0:=\frac{s_1-s_2}{2}$. Then
    \begin{align}\label{eq E(z_n^s1 z_n^s2) growth rate2}
        \E{Z_n^{s_1}(x)Z_n^{s_2}(\overline{x})}
        \sim
        n^{\binom{s_1+s_2}{k_0-1}} \Bigl(C(k_0-1)(x^{k_0-1})^n +
        C(k_0+1)(x^{k_0+1})^n \Bigr).
    \end{align}
\end{itemize}

Additionally, for every even $k_0$ with $-s_2\leq k_0\leq s_1$
$$C(s_1,s_2,-k_0)=\overline{C(s_2,s_1,k_0)}.$$
\end{theorem}
\begin{proof}
We prove here only the case $s_1+s_2=4p+1$ and $s_2$ even or odd. The other cases are similar.\\
We have to calculate
$$M:=\max_{k\text{ even}}\binom{s_1+s_2}{s_2+k}.$$
We know that $\max\binom{4p+1}{k}=\binom{4p+1}{2p}=\binom{4p+1}{2p+1}$. If $s_2$ is even then $k+s_2$ runs through all even
numbers between $0$ and $s_1+s_2$. Therefore $M=\binom{4p+1}{2p}$ and the maximum is attained
for $k_0=2p-s_2=\frac{s_1+s_2-1}{2}-s_2=\frac{s_1-s_2-1}{2}$. We have in this case $s_1-s_2=4p+1-2s_2=4m+1$ and
formula \eqref{eq E(z_n^s1 z_n^s2) growth rate1} follows from \eqref{growth_rate}. The argument for $s_2$ odd is similar.\\
It remains to show that
$$
C(s_1,s_2,-k_0)=C(s_2,s_1,k_0).
$$
This follows from
\begin{align*}
C(s_1,s_2,-k_0,x)
&=
\frac{1}{\Bigl(\binom{s_1+s_2}{s_2-k_0}-1\Bigr)!}
\prod_{\substack{k=-s_2\\k\neq -k_0}}^{s_1}(1-x^k \overline{x}^{-k_0})^{(-1)^{k+1}\binom{s_1+s_2}{s_2+k}}\\
&=
\frac{1}{\Bigl(\binom{s_1+s_2}{s_1+k_0}-1\Bigr)!}
\prod_{\substack{k=-s_1\\k\neq k_0}}^{s_2}(1-x^{-k} x^{k_0})^{(-1)^{k+1}\binom{s_1+s_2}{s_1+k}}\\
&=
\overline{C(s_2,s_1,k_0,x)}
.
\end{align*}

\end{proof}
%
%
%
%
%
%
%
%
%
%
%
%
%
\begin{corollary}
\label{growth s1=s2}
$$
\E{|Z_n(x)|^{2s}}
\sim
n^{\binom{2s}{s}-1}\ \frac{\prod_{k=1}^s |1-x|^{2\binom{2s}{s+k}}}{\Bigl(\binom{2s}{s}-1\Bigr)!}.
$$
\end{corollary}
\begin{proof}
Put $s_1=s_2=s$ in theorem \ref{thm_z^s1 z^s2 groth}.
\end{proof}
\begin{corollary}
$$
\mathrm{Var}\Bigl(Z_n(x)\Bigr)
\sim
n \ |1-x|^2.
$$
\end{corollary}
\begin{proof}
We have $\E{|Z_n(x)|^{2}}\sim C(1,1,0,x) n$ and $\E{Z_n(x)}=1-x$ (see \eqref{exp z_n}).
\end{proof}
%
%
%
%
%
%
\subsection{The real and the imaginary part}
\label{last_result}
We mentioned in the introduction the results in \cite{HKOS}. Do we have the same results for $Z_n(x)$?\\
We first look at the expectation and the variance of the real and the imaginary of $Z_n(x)$.
We set $R_n(x):=\Re\Bigr(Z_n(x)\Bigl)$ and $I_n(x):=\Im\Bigr(Z_n(x)\Bigl)$. We have
\begin{lemma}
\label{real variance}We write $x=e^{i\varphi}$. Then
\begin{enumerate}
  \item $\label{normal 1}\E{R_n(x)}=1-\cos(\varphi)$,
  \item $\label{normal 2}\E{I_n(x)}=-\sin(\varphi)$,
  \item $\label{normal 3}\mathrm{Var}\Bigl(R_n(x)\Bigr) \sim n \frac{|1-x|^2}{2}$,
  \item $\label{normal 4}\mathrm{Var}\Bigl(I_n(x)\Bigr) \sim n \frac{|1-x|^2}{2}$,
  \item $\label{normal 5}\mathrm{Corr}(R_n,I_n)\rw 0$ for $n\rw\infty$.
\end{enumerate}
\end{lemma}
\begin{proof}
\eqref{normal 1} and \eqref{normal 2} follows from \eqref{exp z_n}.\\
As next we prove \eqref{normal 3} and \eqref{normal 4}.
We use the growth rates for $s_1=s_2=1$ and $s_1=2, s_2=0$. We only give the important constants explicitly.
\begin{subequations}
     \begin{align}
       \E{Z_n(x)Z_n(\overline{x})}=&\E{R^2_n+I_n^2}                &\sim  |1-x|^2 n         \label{calc var 1} &,\\
       \E{Z^2_n(x)}               =&\E{R^2_n+2i R_nI_n-I_n^2}      &\sim  C_1 + C_2 (x^2)^n \label{calc var 2} &.
     \end{align}
\end{subequations}
We calculate $\eqref{calc var 1}\pm\Re\eqref{calc var 2}$ and get
$$\E{2R_n^2}\sim |1-x|^2 n$$
$$\E{2I_n^2}\sim |1-x|^2 n$$
We now prove the last point. We know from \eqref{calc var 2} that
$$\mathrm{Cov}(R_n,I_n)
=\E{R_n I_n}+\sin(\varphi)(1-\cos(\varphi))
\sim C_4 + C_5 \sin(2n\varphi) + C_6 \cos\left(2n\varphi\right).$$
The point \eqref{normal 5} now follows from \eqref{normal 3} and \eqref{normal 4}.
\end{proof}
What are the growth rates or $\E{R_n^s}$ and $\E{I_n^s}$? We need the following lemma to answer this question.
\begin{lemma}
\label{lem_r^s=sum z^k z^s-k}
Let $s\in\Nr$ and $z=x+iy$ be given with $x,y\in \Rr$. Then
\begin{align}
x^s
=
\frac{1}{2^s}\sum_{k=0}^s \binom{s}{k} z^k \overline{z}^{s-k};
\qquad
y^s
=
\frac{1}{(2i)^s} \sum_{k=0}^s (-1)^{s+k}\binom{s}{k} z^k \overline{z}^{s-k}
\end{align}
\end{lemma}
\begin{proof}
We argue with induction. If $s=1$ then $x=\frac{1}{2}(z+\overline{z})$.\\
$s\rw s+1$: We have
$$x^{s+1}
=
x\ x^s
=
\frac{1}{2}(z+\overline{z})\frac{1}{2^s}\sum_{k=0}^s \binom{s}{k} z^k \overline{z}^{s-k}
=
\frac{1}{2^{s+1}}\sum_{k=0}^{s+1} \binom{s+1}{k} z^k \overline{z}^{(s+1)-k}.$$
The proof for $y^s$ is similar.
\end{proof}

%
%
%
We then have
\begin{theorem}
\label{thm_growth Re and IM x=1}
Choose any $s\in \Nr$ and write $x=e^{i\varphi}$ with $\varphi\in[0,2\pi]\setminus 2\pi\Q$.
Then there exists (real)
constants $a_{2k}=a_{2k}(\varphi,s)$ and $b_{2k}=b_{2k}(\varphi,s)$ for $0\leq k\leq \lfloor (s+1)/4\rfloor$ with
\begin{align}\label{eq_growth_R^s}
\E{R_n^s}
&\sim n^{\binom{s}{\lfloor s/2\rfloor}} \left(\sum_{k=0}^{\lfloor (s+1)/4\rfloor} a_{2k} \cos\bigl((2k)n\varphi\bigr)\right),\\
\E{I_n^s}
&\sim n^{\binom{s}{\lfloor s/2\rfloor}} \left(\sum_{k=0}^{\lfloor (s+1)/4\rfloor} b_{2k} \sin\bigl((2k)n\varphi\bigr)\right).
\end{align}
At least one $a_{2k}$ and one $b_{2k}$ is not equal to zero.
\end{theorem}
\begin{proof}
We only prove the behavior for $\E{R_n^s}$ and $s=4m$. The other cases are similar.\\
We have
\begin{align*}
\E{R_n^s}
&=
\frac{1}{2^s} \sum_{k=0}^s\binom{s}{k} \E{Z_n^k(x)Z_n^{s-k}(\overline{x})}\\
&=
\frac{1}{2^s} \left(\sum_{\substack{k=0\\ k\text{ even}}}^s\binom{s}{k}\E{Z_n^k(x)Z_n^{s-k}(\overline{x})}  +\sum_{\substack{k=0\\ k\text{ odd}}}^s\binom{s}{k} \E{Z_n^k(x)Z_n^{s-k}(\overline{x})}\right).
\end{align*}
We now apply theorem \ref{thm_z^s1 z^s2 groth}. If $k$ is odd then $k-(4m-k)=4p+2$ for a $p\in\Zr$ and
the growth rate of $\E{Z_n^k(x)Z_n^{s-k}(\overline{x})}$ is $\binom{4m}{2m+1}(\dots)$. If $k$ is even then $k-(4m-k)=4p$ for a $p\in\Zr$ and
the growth rate of $\E{Z_n^k(x)Z_n^{s-k}(\overline{x})}$ is $\binom{4m}{2m}(\dots)$. It is therefore sufficient to look at even $k$. We get
\begin{align*}
\E{R_n^s}
&\sim
n^{\binom{4m}{2m}-1} \left(\frac{1}{2^s} \sum_{\substack{k=0\\ k\text{ even}}}^s \binom{s}{k} \left(x^{\frac{k-(s-k)}{2}}\right)^n C\left(\frac{k-(s-k)}{2}\right) \right)\\
&\sim
n^{\binom{4m}{2m}-1} \left(\frac{1}{2^s} \sum_{k=-m}^{m} \binom{s}{2m+2k} (x^{2k})^n C\left(2k\right) \right)  \label{eq_bew_groth r^s}\\
&\sim
n^{\binom{4m}{2m}-1} \frac{1}{2^s}
\left(\binom{s}{2m} C(0)+2 \sum_{k=1}^{m} \binom{s}{2m+2k} \Bigl(\Re(C(2k)) \cos\bigl((2k)n\varphi\bigr) -\Im(C(2k)) \sin\bigl((2k)n\varphi\bigr)\Bigr)\right).
\end{align*}
We have used in the last inequality that $C(s_1,s_2,-k_0) = \overline{C(s_2,s_1,k_0)}$.\\

This proves \eqref{eq_growth_R^s} if we can show that the last bracket is equal zero only for finitely many $n$ (and $x$ fixed). \\
We define
$$g(t):=\binom{s}{2m} C(0)+2 \sum_{k=1}^{m} \binom{s}{2m+2k}\Bigl(\Re(C(2k)) \cos\bigl((2k)t\bigr) -\Im(C(2k)) \sin\bigl((2k)t\bigr)\Bigr).$$
Suppose there are infinitely many (different) $n\in\Nr$ with $g(n\varphi)=0$.
All numbers $n\varphi$ are different modulo $2\pi$, since $\varphi\notin 2\pi\Q$.
Therefore there are infinite many $t\in[0,2\pi]$ with $g(t)=0$.
But $g(t)$ is a non trivial linear combination of $\cos(\cdot)$ and $\sin(\cdot)$ and is therefore a holomorphic function in $t$.
It follows immediately from the identity theorem (see \cite{busam}) that $g(t)\equiv 0$.
This is a contradiction since the functions $\cos(m_1t)$ and $\sin(m_2t)$ are linearly independent for $m_1\geq 0,m_2 >0$.
%
\end{proof}
\section{Concluding Remarks}
The behavior of $\E{\left(\frac{d}{dx}Z_n(x)\right)^s}$ for $s\in\Nr,|x|<1$ is not included in this paper.
We have been able to prove a result similar to theorem~\ref{conver x<1}. Our proof uses the techniques of section~\ref{holoins} and Hartog's theorem. Unfortunately we could not give an explicit expression for the limit and the prove is more difficult than the proof for $\E{Z_n^s(x)}$. We therefore decided to omit this result.\\
We have proven in this paper several result, but there are still open questions:
\begin{enumerate}
\item Is $Z_\infty(x)$ uniquely defined by its moments?
\item Are there random variables $R_\infty$ and $I_\infty$ such that $R_n \xrightarrow{d} R_\infty$ and $I_n \xrightarrow{d} I_\infty$?
\item What are the growth rates of $\E{Z_n^{s}(x)}$ for $s\in\C$, $|x|=1$ and $x$ not a root of unity?
\end{enumerate}

\textbf{Acknowledgement}
I would like to thank Ashkan Nikeghbali and Paul-Olivier Dehaye for their help writing my first paper.

\bibliographystyle{plain}

\end{document}